%% file: CubesArxiv.tex
\documentclass{article}

\usepackage{amsmath}
\usepackage{amsthm}
\usepackage{amssymb}
\usepackage{amsfonts}
\usepackage{latexsym}
\usepackage{delarray}
\usepackage[dvips]{graphics}
\usepackage{epsfig}
\usepackage{color}

\newtheorem{thm}{Theorem}[section]
\newtheorem{lemma}[thm]{Lemma}
\newtheorem{prop}[thm]{Proposition}
\newtheorem{cor}[thm]{Corollary}
\newtheorem{question}{Question}

\theoremstyle{definition}

\theoremstyle{remark}

\newcommand{\Aut}{\operatorname{Aut}}
\newcommand{\Stab}{\operatorname{Stab}}

\title{The distinguishing number of the augmented cube and hypercube powers}
\author{Melody Chan \\ University of Cambridge \\ Cambridge, England \\ \texttt{melody.chan@aya.yale.edu}}

\begin{document}
\maketitle

\begin{abstract}
The distinguishing number of a graph $G$, denoted $D(G)$, is the minimum number of colors such that there exists a coloring of the vertices of $G$ where no nontrivial graph automorphism is color-preserving.  In this paper, we answer an open question posed in \cite{b_c} by showing that the distinguishing number of $Q_n^p$, the $p^{\textrm{th}}$ graph power of the $n$-dimensional hypercube, is 2 whenever $2 < p < n-1$.  This completes the study of the distinguishing number of hypercube powers.  We also compute the distinguishing number of the augmented cube $AQ_n$, a variant of the hypercube introduced in \cite{c_s}.  We show that $D(AQ_1) = 2$; $D(AQ_2) = 4$; $D(AQ_3) = 3$; and $D(AQ_n) = 2$ for $n \ge 4$.  The sequence of distinguishing numbers $\{D(AQ_n)\}_{n=1}^{\infty}$ answers a question raised in \cite{a_c_intro}. 
\end{abstract}

\section{Introduction}

Given a graph $G$, an \emph{$r$-coloring} of $G$ is a function $c \colon V(G) \rightarrow \{1, \ldots, r\}$.  An automorphism $\pi$ of the graph $G$ is said to \emph{preserve} the coloring $c$ if $c(\pi v) = c(v)$ for each vertex $v \in V(G)$.  A coloring of $G$ is said to be \emph{distinguishing} if no nontrivial automorphism of $G$ preserves it.  The \emph{distinguishing number} of $G$, denoted $D(G)$, is the smallest number of colors $r$ such that there exists a distinguishing $r$-coloring of $G$.  Note that throughout this paper, $r$-colorings are not required to be proper graph colorings; two adjacent vertices may or may not have the same color.

The distinguishing number was first introduced in \cite{a_c} by Albertson and Collins, who proved some general results relating the distinguishing number of a graph to properties of its automorphism group.  For example, they showed that $D(G) \le 2$ if $\Aut(G)$ is abelian, and $D(G) \le 3$ if $\Aut(G)$ is dihedral.  They also proved that $D(G) = 2$ or $D(G) = 4$ if $\Aut(G) \cong S_4$.  In \cite{cctcheng}, Cheng gave an algorithm for computing the distinguishing number of trees and forests.  Potanka computed the distinguishing number of the generalized Petersen graphs in \cite{potanka}.  In addition, Russell and Sundaram considered the computational complexity of the distinguishing number in \cite{r_s}.

In \cite{b_c}, Bogstad and Cowen computed $D(Q_n)$ and $D(Q_n^2)$ for each $n$, where $Q_n$ denotes the $n$-dimensional hypercube and $Q_n^p$ denotes its $p^{\textrm{th}}$ graph power.  They also noted that $D(Q_n^p)$ is easily computable when $p \ge n-1$.  They left $D(Q_n^p)$ for $2 < p < n-1$ as an open question, conjecturing that $D(Q_n^p) = 2$ for fixed $p$ and sufficiently large $n$.  Additionally, they offered the weaker conjecture that there exists a constant $c$ such that $D(Q_n^p) \le cp$ for fixed $p$ and sufficiently large $n$.

In this paper, we show that $D(Q_n^p) = 2$ for each $p$ and $n$ satisfying $2 < p < n-1$.  The proof relies on a surprising theorem in \cite{m_p} on the automorphism group of hypercube powers.  This result completes the determination of the distinguishing number of all hypercube powers of all dimensions.

We then move on to consider the distinguishing number of the augmented cube, introduced by Choudum and Sunitha in \cite{c_s} as a variation on the hypercube possessing several favorable network properties.  We compute the distinguishing number of the augmented cube of each dimension.  Finally, we answer an open question posed by Albertson and Collins in \cite{a_c_intro} on the existence of a class of graphs $\{G_n\}_{n=1}^{\infty}$ whose sequence of distinguishing numbers increases to some $k \ge 4$ and then decreases to 2.  We show that the augmented cubes have precisely this property.

\section{Hypercubes} \label{qn}

The $n$-dimensional hypercube, denoted $Q_n$, is the graph on $2^n$ vertices labeled by length-$n$ binary strings $\{a_1 \cdots a_n~|~a_i \in \{0,1\}\}$ and where two vertices $v = v_1 \cdots v_n$ and $w = w_1 \cdots w_n$ are joined by an edge if and only if $|\{ i~|~v_i \ne w_i\}| = 1$.  The $p^{th}$ power of a graph $G$, denoted $G^p$, is defined to be a new graph with the same vertex set as $G$ and in which two distinct vertices in $G^p$ are connected by an edge if the corresponding vertices in $G$ are at distance at most $p$.  More formally, $V(G^p) = V(G)$ and $E(G^p) = \{ \{v,w\}~|~0 < d_G(v,w) \le p \}$.  Here $d_G(v,w)$ denotes the length of the shortest path between $v$ and $w$ in $G$.  Thus $G^1 \cong G$, and if $p$ is at least the diameter of $G$ then $G^p \cong K_{|V(G)|}$, the complete graph on $|V(G)|$ vertices.

In \cite{b_c}, Bogstad and Cowen consider the distinguishing number of the hypercube and the second power of the hypercube.  For $n \in \{2,3\}$, they prove that $D(Q_n) = 3$ and $D(Q_n^2) = 4$.  For $n \ge 4$, they show $D(Q_n) = D(Q_n^2) = 2$.   They note further that the graph $Q_n^{n-1}$ consists of the complement of a perfect matching on $2^n$ vertices, and both graphs have distinguishing number $\min\{x~|~\binom{x}{2} \ge 2^{n-1} \}$.  Finally, we have already seen that for $p \ge n$, the graph $Q_n^p$ is isomorphic to the complete graph $K_{2^n}$ and so has distinguishing number $2^n$.

The authors leave $D(Q_n^p)$ for $2 < p < n-1$ as an open question.  At this point, we wish to draw the reader's attention to the following very surprising theorem proved in \cite{m_p}.

\begin{thm} \cite[Section 1]{m_p} \label{t:m_p} For $2 < p < n-1$, 
$$
\Aut(Q_n^p) = 
\begin{cases}
\Aut(Q_n) & \textrm{if $p$ is odd,} \\
\Aut(Q_n^2) & \textrm{if $p$ is even}.
\end{cases}
$$
\end{thm}  

For clarity's sake, we note the following subtlety.  It is shown in \cite{a_c} that two graphs with automorphism groups that are isomorphic may still have different distinguishing numbers.  However, Theorem \ref{t:m_p} gives more than just isomorphisms between the groups under consideration.  For note that an automorphism of $Q_n$ preserves all distances in $Q_n$ and is therefore an automorphism of $Q_n^p$ for any power $p$.  So $\Aut(Q_n)$ is realized as a subgroup of $\Aut(Q_n^p)$.  Theorem \ref{t:m_p} tells us that for any odd $p$, $\Aut(Q_n)$ and $\Aut(Q_n^p)$ are in fact precisely the same subgroup of the permutation group of their vertices, and so they act with equal distinguishing number.   Thus $D(Q_n^p) = D(Q_n) = 2$ for odd $p$.  A similar argument shows that $D(Q_n^{2k}) = D(Q_n^2) = 2$.  We summarize as follows:

\begin{cor}
$D(Q_n^p) = 2$ for $2 < p < n-1$.
\end{cor}
This gives a complete answer to the question posed in \cite{b_c}.

Before concluding our discussion of hypercube powers, we state a simple but useful lemma.

\begin{lemma} \label{l:subgroup_graph}
Suppose $G_1$ and $G_2$ are graphs on the same vertex set, and $\Aut(G_1)$ is a subgroup of $\Aut(G_2)$.  Then $D(G_1) \le D(G_2)$.
\end{lemma}

\begin{proof} By definition, there exists a $D(G_2)$-coloring of the vertices of $G_2$ such that no nonidentity automorphism of $G_2$ preserves it.  In particular, since $\Aut(G_1) \le \Aut(G_2)$, no nonidentity automorphism of $G_1$ preserves this coloring.  So $D(G_2)$ colors suffice to produce a distinguishing coloring of $G_1$.
\end{proof}

Lemma \ref{l:subgroup_graph} shows that one of the main theorems in \cite{b_c} implies another.  Indeed, it is shown that $D(Q_n^2) = 2$ for $n \ge 4$.  However, it is clear that $\Aut(Q_n)$ is a subgroup of $Aut(Q_n^2)$, so $D(Q_n) \le 2$.  Since $\Aut(Q_n)$ is nontrivial, $D(Q_n) > 1$, so $D(Q_n) = 2$ for $n \ge 4$.

\section{Augmented cubes} \label{aqn}

The $n$-dimensional augmented cube, denoted $AQ_n$, is a hypercube variant introduced in \cite{c_s} by Choudum and Sunitha.  As with the hypercube, the vertices of $AQ_n$ are length-$n$ binary strings $\{a_1 \cdots a_n~|~a_i \in \{0,1\}\}$.  The edges of the augmented $n$-cube, however, are a superset of the edges of the $n$-cube.  We define $AQ_n$ recursively as follows.  For $n=1$, let $AQ_1 \cong K_2$.  To construct $AQ_n$ for $n>1$, we take two copies of $AQ_{n-1}$ and connect not only pairs of corresponding vertices, as in the hypercube, but also pairs of opposite vertices.  More precisely, let us index our copies of $AQ_{n-1}$ as $AQ_{n-1}^0$ and $AQ_{n-1}^1$, with vertex sets $V(AQ_{n-1}^0) = \{0a_2 \cdots a_n~|~a_i \in \{0,1\}\}$ and $V(AQ_{n-1}^1) = \{1b_2 \cdots b_n~|~b_i \in \{0,1\}\}$.  We add an edge between vertices $a = 0a_2 \cdots a_n \in AQ_{n-1}^0$ and $b = 1b_2 \cdots b_n \in AQ_{n-1}^1$ if either 

(1) $a_i = b_i$ for each $2 \le i \le n$, or 
\nopagebreak

(2) $a_i \ne b_i$ for each $2 \le i \le n$.  

Thus, $AQ_2$ is isomorphic to $K_4$, the complete graph on 4 vertices.  The augmented 3-cube, $AQ_3$, is shown in Figure \ref{aq3}.  We note that $AQ_n$ is a $(2n-1)$-regular graph with diameter $\lceil \frac{n}{2} \rceil$.  

A useful characterization of adjacency that follows directly from the recursive definition of $AQ_n$ is as follows.

\begin{prop} \cite[Proposition 2.1]{c_s} \label{p:aqnadj}
The vertices $a = a_1 \cdots a_n$ and $b = b_1 \cdots b_n$ are adjacent in $AQ_n$ if and only if

(1) there exists $l$, $1 \le l \le n$, such that $a_i = b_i$ for $i \ne l$ and $a_l \ne b_l$, or

(2) there exists $l$, $1 \le l \le n$, such that for $1 \le i \le l-1$, $a_i = b_i$, and for $l \le i \le n$, $a_i \ne b_i$.
\end{prop}

In what follows, we compute $D(AQ_n)$ for each $n$.  First, we present a lemma that is true for each $n$ but will be used when $n \ge 3$.  Throughout, we let $V = V(AQ_n)$ and $E = E(AQ_n)$.  Also, we let $\bar{x}_i = 1 - x_i$ for $x_i \in \{0, 1\}$.  Finally, we denote the vectors $0 \cdots 0$ and $0 \cdots 0 1$ by $\mathbf{0}$ and $\mathbf{1}$ respectively.

\begin{lemma} \label{l:mainlemma}  
Fix $n$ and suppose a coloring $c$ of $AQ_n$ has the property that for any two vertices $x = x_1 \cdots x_n$ and $y = y_1 \cdots y_n$ satisfying $x_n \ne y_n$ and both different from $\mathbf{0}$ and $\mathbf{1}$, we have $c(x) \ne c(y)$.  Suppose further that a graph automorphism $\pi$ is color-preserving with respect to $c$ and fixes $\mathbf{0}$ and $\mathbf{1}$.  Then $\pi$ is the identity automorphism.
\end{lemma}

\begin{proof} 

For $1 \le i \le n$, let $B_i$ be the subgraph induced by the vertices in the set \linebreak[4] $\{0 \cdots 0~x_{n-i+1} \cdots x_n~|~x_i \in \{ 0 , 1 \} \}$.  Thus, for each $i$, 
$B_i \cong AQ_i$ and $B_1 \subset B_2 \subset \cdots \subset B_n = AQ_n$.  
We will prove by induction on $i$ that $\pi$ fixes each vertex of $B_i$ for $1 \le i \le n$.  
The case $i=1$ is true by assumption.  
Now suppose $\pi$ fixes each vertex of $B_i$.  We wish to show that $\pi$ also fixes each $v = v_1 \cdots v_n \in B_{i+1} \setminus B_i$.  Here, $v$ must have the form $v = 0 \cdots 0 1 v_{n-i+1} \cdots v_n$.  
Thus, $v$ differs from every vertex in $B_i$ in coordinate $n-i$.  Then by Proposition~\ref{p:aqnadj}, $v$ is adjacent to precisely two vertices in $B_i$, namely $\alpha = 0 \cdots 0 v_{n-i + 1} \cdots v_n$, where $v$ and $\alpha$ differ in coordinate $n-i$ only, and $\beta = 0 \cdots 0 \bar{v}_{n-i+1} \cdots \bar{v}_n$, where $v$ and $\beta$ differ in coordinate $n-i$ and every subsequent coordinate.  We claim that $v$ is the sole vertex in $V \setminus B_i$ of color $c(v)$ and adjacent to both $\alpha$ and $\beta$, and therefore that $v$ must be fixed by $\pi$.

First, consider the other vertices in $B_{i+1} \setminus B_i$.  By Proposition~\ref{p:aqnadj}, there is precisely one other vertex in $B_{i+1} \setminus B_i$ adjacent to both $\alpha$ and $\beta$, namely $v' = 0 \cdots 0 1 \bar{v}_{n-i+1} \cdots \bar{v}_n$.  But $v$ and $v'$ differ in their last coordinate (and neither equals $\mathbf{0}$ or $\mathbf{1}$) so that $c(v) \ne c(v')$ by assumption.  

Next, consider the vertices in $V \setminus B_{i+1}$.  We claim that none of these vertices is adjacent to both $\alpha$ and $\beta$.  Suppose for a contradiction that there exists $w = w_1 \cdots w_n \in V \setminus B_{i+1}$ adjacent to $\alpha$ and $\beta$.  Since $w \not\in B_{i+1}$, we have $w_j = 1$ for some $j < n-i$ and thus $w_j \ne {\alpha}_j$ and $w_j \ne {\beta}_j$.  Now, $w_n$ differs from one of ${\alpha}_n$ and ${\beta}_n$ since ${\alpha}_n \ne {\beta}_n$, so without loss of generality assume $w_n \ne {\alpha}_n$.  Then by Proposition~\ref{p:aqnadj}, since $w$ and $\alpha$ differ in coordinates $j$ and $n$, they must also differ in every coordinate between $j$ and $n$; in particular $w_{j+1} \ne {\alpha}_{j+1}$.  But ${\alpha}_{j+1} = {\beta}_{j+1} = 0$ since $j+1 \le n-i$.  So $w_{j+1} \ne {\beta}_{j+1}$.  Then $w$ and $\beta$ differ in coordinates $j$ and $j+1$, so by Proposition~\ref{p:aqnadj}, they must also differ in every subsequent coordinate and in particular in coordinate $n$.  Thus $w_n \ne {\beta}_n$.  But $w_n \ne {\alpha}_n$, ${\alpha}_n \ne {\beta}_n$ and all three are in $ \{ 0 , 1 \} $ so we have a contradiction.  

Therefore, $v$ is the only vertex in $V \setminus B_i$ of color $c(v)$ and adjacent to both $\alpha$ and $\beta$ in $B_i$.  Since $\pi$ fixes each vertex of $B_i$, $\pi$ must fix $v$ as well.  Thus, every vertex of $B_{i+1}$ is fixed under $\pi$.  Finally, we proceed by induction to conclude that $\pi$ must fix every vertex of $B_n = AQ_n$ and therefore that $\pi$ is the identity automorphism.
\end{proof}

Now we are ready to state the main theorem of the section.

\begin{thm} \label{t:aqn_mainthm}
$$
D(AQ_n) =
\begin{cases}
2 & \textrm {if~~}  n=1 \\
4 & \textrm{if~~}n=2 \\
3 & \textrm{if~~} n=3 \\
2 & \textrm{if~~} n \ge 4
\end{cases}
$$
\end{thm}

\begin{proof}  The cases $n=1$ and $n=2$ follow immediately from the fact that $AQ_1 \cong K_2$ and $AQ_2 \cong K_4$, and we have $D(K_n) = n$ for all $n$.  We will now consider the cases $n=3$ and $n \ge 4$ separately.

\begin{lemma} $D(AQ_3) = 3$.
\end{lemma}

\begin{figure}
\begin{center}
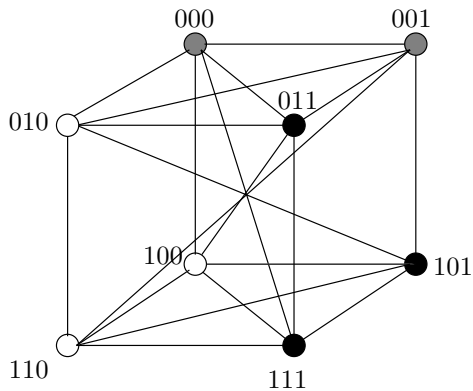
\caption{A distinguishing 3-coloring of $AQ_3$.}
\label{aq3}
\end{center}
\end{figure}

\begin{proof}
We will exhibit a distinguishing 3-coloring of $AQ_3$, and then show that no distinguishing 2-coloring exists.  Let $c\colon  V(AQ_3) \rightarrow \{ 1, 2, 3 \}$ be the coloring $c(000) = c(001) = 3$, $c(010) = c(100) = c(110) = 1$, $c(011) = c(101) = c(111) = 2$.  This coloring is shown in Figure \ref{aq3}, where colors 1, 2, and 3 correspond to white, black, and gray, respectively.  We claim that $c$ is a distinguishing 3-coloring of $AQ_3$.  First we show that a color-preserving graph automorphism $\pi$ must fix $000$ and $001$ pointwise.  Consider the subgraph of $AQ_3$ induced by the vertices of colors $1$ and $3$.  Of course the restriction of $\pi$ to this subgraph must also be an automorphism of it.  Note that vertex $100$ is the unique vertex of degree 2 in this subgraph.  Therefore, $\pi$ must fix it.  Since $100$ is adjacent to $000$ but not $001$ in $AQ_3$, and $000$ and $001$ are the only vertices of color 3, $\pi$ must fix them as well.  Finally, we apply Lemma~\ref{l:mainlemma} to conclude that $\pi$ must be the identity automorphism and therefore that $D(AQ_3) \le 3$.

It remains to be shown that $D(AQ_3) > 2$.  Suppose instead that $c\colon  V(AQ_3) \rightarrow \{1, 2 \}$ is a distinguishing 2-coloring of $V(AQ_3)$.  We will produce a contradiction by constructing a non-trivial automorphism of $AQ_3$ that preserves $c$.  Let $\mathbf{0}$ denote the vertex $000$, and for each vertex $x = x_1 x_2 x_3$, let $x^*$ denote the vertex $x_1 \bar{x}_2 \bar{x}_3$.  The main observation in this proof (one that does not generalize nicely to higher dimensions) is that each $x$ is adjacent to precisely the same set of vertices as $x^*$.  This observation can be checked case by case using Proposition~\ref{p:aqnadj}.  Thus the transposition $(x \textrm{ } x^*)$, for each pair $x$ and $x^*$, is an automorphism of $AQ_3$.  Then we must have $c(x) \ne c(x^*)$ for each $x$.  Thus, there are 4 vertices of color 1 and 4 vertices of color 2.  Without loss of generality, suppose $c(\mathbf{0}) = 1$, then $c(\mathbf{0}^*) = 2$.  Let $v \not\in \{\mathbf{0}, \mathbf{0^*}\}$ be another vertex colored 1, then $c(v^*) = 2$.  
Let $\pi\colon  V(AQ_3) \rightarrow V(AQ_3)$ be given by $\pi(x) = x + v$, where addition is carried out coordinate-wise in $\mathbb{F}_2$.  By Proposition \ref{p:aqnadj}, $\pi \in \Aut(AQ_3)$.
We may express $\pi$ in cyclic notation as $(\mathbf{0} \textrm{ } v) (\mathbf{0}^* \textrm{ } v^*) (a \textrm{ } b) (d \textrm{ } e)$, where $\{ a, b, d, e \} = V(AQ_3) \setminus \{ \mathbf{0}, \mathbf{0}^*, v, v^* \}$.  
 Now, we know that two vertices of $\{ a, b, d, e \}$ are colored 1 and two vertices are colored 2.  If $c(a) = c(b)$ then $c(d) = c(e)$, and $\pi$ is a non-trivial color preserving automorphism, which contradicts that $c$ is a distinguishing coloring.  Thus we may assume, without loss of generality, that $c(a) = c(d) = 1$ and $c(b) = c(e) = 2$.  Now, since $c(a^*) = 2$, either $b = a^*$ or $e = a^*$.  But $v \ne \mathbf{0}^*$, so $b = \pi(a) = a + v \ne a + \mathbf{0}^* = a^*$.  So $e = a^*$.  Likewise, $d = b^*$.  Then let $\tau_1 = (a \textrm{ } e)$, $\tau_2 = (b \textrm{ } d)$, both of which are in $\Aut(AQ_3)$.  Then the composition $\tau_1 \tau_2 \pi = (\mathbf{0} \textrm{ } v) (\mathbf{0}^* \textrm{ } v^*) (a \textrm{ } d) (b \textrm{ } e)$ is a nontrivial automorphism of $AQ_3$ preserving $c$, contradicting the assumption that $c$ is a distinguishing coloring.  Thus $D(AQ_3) > 2$.  We have already established that $D(AQ_3) \le 3$, so therefore $D(AQ_3) = 3$.
\end{proof}

\begin{lemma} $D(AQ_n) = 2$ for $n \ge 4$.
\end{lemma}

\begin{proof}
Let $c\colon  V(AQ_n) \rightarrow \{1, 2 \}$ be given by
$c(\mathbf{0}) = 2$, $c(\mathbf{1}) = 1$, and for $x = x_1 \cdots x_n$ different from $\mathbf{0}$ and $\mathbf{1}$, 
$c(x_1 \cdots x_n) = x_n + 1$.  We claim that $c$ is a distinguishing 2-coloring of $AQ_n$.
We will show that any $\pi$ that preserves our coloring $c$ fixes $\mathbf{0}$ and $\mathbf{1}$, and apply Lemma~\ref{l:mainlemma} to conclude that $\pi$ must be the identity automorphism.  Then since $AQ_n$ has non-trivial automorphism group as given in \cite{c_stechnical}, we have that $D(AQ_n) = 2$.

Associate with each vertex $w$ the ordered triple of natural numbers $t_w = (x_1, x_2, x_3)$ where $x_1 = c(w)$ is the color of $w$, where $x_2 = |\{v \in V~|~(v,w) \in E \textrm{ and } c(v) = 1 \}|$ is the number of color-1 neighbors of $w$, and where $x_3 = |\{v \in V~|~(v,w) \in E \textrm{ and } c(v) = 2\}|$ is the number of color-2 neighbors of $w$.  For $i = \mathbf{0},\mathbf{1}$ and $j = 1,2$, let 
$M_{i,j}$ be the set $\{v \in V~|~v \in N[i] \textrm{ and } c(v) = j \}$.  Here, $N[i] = \{v \in V~|~(v, i) \in E \} \cup \{ i \}$ denotes the closed neighborhood of a vertex $i$.  Note that a given vertex $w = w_1 \cdots w_n \in V$ has $n-1$ neighbors with last coordinate $w_n$, namely
$\bar{w}_1 w_2 \cdots w_n$, $w_1 \bar{w}_2 w_3 \cdots w_n$, $\ldots$, and $w_1 \cdots w_{n-2} \bar{w}_{n-1} w_n$;
and $n$ neighbors with last coordinate $\bar{w}_n$, namely
$w_1 \cdots w_{n-1} \bar{w}_n$, $w_1 \cdots w_{n-2} \bar{w}_{n-1} \bar{w}_{n}$, $\ldots$, and $\bar{w}_1 \cdots \bar{w}_n$.  Using this fact, the following equations are straightforward to check:

\begin{align}
t_{w} & = (1, n-2, n+1) &&\textrm{ for } w \in M_{\mathbf{0},1} \setminus M_{\mathbf{1},1} \\
t_{w} & = (2, n-1, n) &&\textrm{ for } w \in M_{\mathbf{0},2} \setminus M_{\mathbf{1},2} \\
t_{w} & = (1, n, n-1) &&\textrm{ for } w \in M_{\mathbf{1},1} \setminus M_{\mathbf{0},1} \\
t_{w} & = (2, n+1, n-2) &&\textrm{ for } w \in M_{\mathbf{1},2} \setminus M_{\mathbf{0},2} \\
t_{w} & = (1, n-1, n) &&\textrm{ for remaining } w \textrm{ of color 1} \label{eqn5}\\
t_{w} & = (2, n, n-1) &&\textrm{ for remaining } w \textrm{ of color 2} \label{eqn6}
\end{align}

It should be noted that Equation~\eqref{eqn5} comprises the cases that $w$ is of color 1 and (i) $w = \mathbf{1}$, (ii) $w$ is adjacent to neither $\mathbf{0}$ nor $\mathbf{1}$, or (iii) $w = 0 \cdots 010$ is adjacent to both $\mathbf{0}$ and $\mathbf{1}$.  Likewise, Equation~\eqref{eqn6} comprises the cases that $w$ is of color 2 and (i) $w = \mathbf{0}$, (ii) $w$ is adjacent to neither $\mathbf{0}$ nor $\mathbf{1}$, or (iii) $w = 0 \cdots 011$ is adjacent to both $\mathbf{0}$ and $\mathbf{1}$.

Now we give an argument showing that $\pi$ must fix $\mathbf{0}$.  The argument that $\pi$ must fix $\mathbf{1}$ is similar and will therefore be omitted.
Note that

\begin{eqnarray*}
M_{\mathbf{0},1} \setminus M_{\mathbf{1},1} &=& \{ 0 \cdots 0 1 0 0 , 0 \cdots 0 1 0 0 0, \ldots, 1 0 \cdots 0 \} \\
M_{\mathbf{0},2} \setminus M_{\mathbf{1},2} &=& \{ 0 \cdots 0 1 1 1 , 0 \cdots 0 1 1 1 1, \ldots, 1 \cdots 1 \}.
\end{eqnarray*}

Let $M = (M_{\mathbf{0},1} \setminus M_{\mathbf{1},1}) \cup (M_{\mathbf{0},2} \setminus M_{\mathbf{1},2})$. 
We see that $M = \{ v \in V~|~t_v = (1, n-2, n+1) \textrm{ or } t_v = (2, n-1, n) \}$, and therefore that $\pi$ must fix $M$ as a set.  Now, $\mathbf{0}$ is adjacent to each vertex in $M$.  In particular, it is adjacent to $x = 0 1 0 \cdots 0$, $y = 1 0 \cdots 0$, $z = 0 1 \cdots 1$, and $u = 1 \cdots 1$.  (Here we use the fact that $d \ge 4$ to guarantee $|M| \ge 4$ and that the vectors $x,y,z,u$ are indeed appropriately defined).  Now, we claim $\mathbf{0}$ is the only vertex not in $M$ that is adjacent to every vertex in $M$.   

Indeed, suppose instead that there exists $b = b_1 \cdots b_n \not\in M \cup \{ \mathbf{0} \}$ and $b$ is adjacent to each of $x,y,z,u$.
We make repeated use of Proposition~\ref{p:aqnadj} in the following cases.

\emph{Case 00}: $b = 0 0 b_3 \cdots b_n$.  Since $b$ and $u$ are adjacent and differ in their first 2 coordinates, they must differ in every coordinate.  Therefore $b = 0 \cdots 0 = \mathbf{0}$, a contradiction.

\emph{Case 01}: $b = 0 1 b_3 \cdots b_n$.  Since $b$ and $y$ are adjacent and differ in their first 2 coordinates, they must differ in every coordinate.  Therefore $b = 0 1 1 \cdots 1 = z \in M$, a contradiction. 

\emph{Case 10}: $b = 1 0 b_3 \cdots b_n$.  Since $b$ and $z$ are adjacent and differ in their first 2 coordinates, they must differ in every coordinate.  Therefore $b = 1 0 \cdots 0 = y \in M$, a contradiction.

\emph{Case 11}: $b = 1 1 b_3 \cdots b_n$.  Since $b$ and $z$ are adjacent, differ in their first coordinate, and share their second  coordinate, they must share every subsequent coordinate.  Therefore $b = 1 \cdots 1 = u \in M$, a contradiction.

Therefore, $\mathbf{0}$ is the only vertex not in $M$ adjacent to each $v \in M$.  Since $\pi$ fixes $M$ as a set, it must fix $\mathbf{0}$.  A similar argument (in which the last bit of each vector is flipped and the two colors are permuted) shows that $\pi$ must fix $\mathbf{1}$.  We apply Lemma~\ref{l:mainlemma} to complete the proof that $D(AQ_n) = 2$ for $n \ge 4$.
\end{proof}

This concludes the proof of Theorem~\ref{t:aqn_mainthm}.

\end{proof}

In \cite{a_c_intro}, Albertson and Collins ask whether there exists a class of graphs $\{G_n\}_{n=1}^{\infty}$ such that the sequence of distinguishing numbers $\{D(G_n)\}_{n=1}^{\infty}$ grows to some $k \ge 4$ and then decreases to 2.  The augmented cubes $\{AQ_n\}_{n=1}^{\infty}$ have precisely this property, as shown in Theorem~\ref{t:aqn_mainthm}.

\section{Discussion and open questions}

Hypercubes and augmented cubes are just two of many classes of graphs for which computing the distinguishing number would be of intrinsic interest.  In addition, one could ask questions relating the distinguishing number to specific graph properties.  The following general question appears in \cite{archive}.

\begin{question} Characterize graphs with distinguishing number 2.
\end{question}

In particular, Saks asks whether a graph that has a nontrivial automorphism group containing no involutions must have distinguishing number greater than 2.

In \cite{jt}, Tymoczko generalizes the notion of the distinguishing number to group actions.  Given a group $\Gamma$ acting on a set $X$, we define the distinguishing number of this action, denoted $D_\Gamma(X)$, to be the smallest number of colors admitting a coloring such that the only elements of $\Gamma$ that induce color-preserving permutations of $X$ are those lying in $\Stab(X)$, the element-wise stabilizer of $X$.  Note that in this case, there exists a faithful action of the group $\Gamma/\Stab(X)$ on $X$ with equal distinguishing number, so we may restrict our attention to faithful actions without loss of generality. Tymoczko shows that the problem of distinguishing group actions is a more general one than distinguishing graphs; for example, there exists a faithful action of $S_4$ with distinguishing number 3, whereas Albertson and Collins proved in \cite{a_c} that no graph with automorphism group $S_4$ has distinguishing number 3.  This leads us to ask the following.

\begin{question}
Given a group $\Gamma$, what integers are realized as distinguishing numbers of faithful actions of $\Gamma$ but not as distinguishing numbers of graphs with automorphism group $\Gamma$?
\end{question}

There seem to be many further interesting questions on the distinguishing number of group actions.  We refer the reader to \cite{mc_products} and \cite{mc_maximum}.

\section{Acknowledgments}

This research was conducted at the University of Minnesota Duluth.  The author would like to express gratitude to Melanie Wood and Philip Matchett for numerous suggestions and help with early drafts of this paper, David Moulton for several helpful conversations, and Joseph Gallian for his support.  Funding was provided by the National Science Foundation (DMS-0137611) and the National Security Agency (H-98230-04-1-0050).

\end{document}

%% file: aqn.pstex_t
\begin{picture}(0,0)%
\includegraphics{aqn.eps}%
\end{picture}%
\setlength{\unitlength}{3947sp}%
\begingroup\makeatletter\ifx\SetFigFont\undefined%
\gdef\SetFigFont#1#2#3#4#5{%
  \reset@font\fontsize{#1}{#2pt}%
  \fontfamily{#3}\fontseries{#4}\fontshape{#5}%
  \selectfont}%
\fi\endgroup%
\begin{picture}(2904,2391)(2326,-2390)
\put(3951,-2390){\makebox(0,0)[lb]{\smash{{\SetFigFont{10}{12.0}{\familydefault}{\mddefault}{\updefault}{\color[rgb]{0,0,0}111}%
}}}}
\put(3366,-116){\makebox(0,0)[lb]{\smash{{\SetFigFont{10}{12.0}{\familydefault}{\mddefault}{\updefault}{\color[rgb]{0,0,0}000}%
}}}}
\put(4731,-116){\makebox(0,0)[lb]{\smash{{\SetFigFont{10}{12.0}{\familydefault}{\mddefault}{\updefault}{\color[rgb]{0,0,0}001}%
}}}}
\put(2326,-766){\makebox(0,0)[lb]{\smash{{\SetFigFont{10}{12.0}{\familydefault}{\mddefault}{\updefault}{\color[rgb]{0,0,0}010}%
}}}}
\put(4016,-635){\makebox(0,0)[lb]{\smash{{\SetFigFont{10}{12.0}{\familydefault}{\mddefault}{\updefault}{\color[rgb]{0,0,0}011}%
}}}}
\put(3171,-1610){\makebox(0,0)[lb]{\smash{{\SetFigFont{10}{12.0}{\familydefault}{\mddefault}{\updefault}{\color[rgb]{0,0,0}100}%
}}}}
\put(4991,-1675){\makebox(0,0)[lb]{\smash{{\SetFigFont{10}{12.0}{\familydefault}{\mddefault}{\updefault}{\color[rgb]{0,0,0}101}%
}}}}
\put(2326,-2325){\makebox(0,0)[lb]{\smash{{\SetFigFont{10}{12.0}{\familydefault}{\mddefault}{\updefault}{\color[rgb]{0,0,0}110}%
}}}}
\end{picture}%

%% file: CubesArxiv.bbl
\begin{thebibliography}{9}

\bibitem{a_c_intro} M. Albertson and K. Collins, An introduction to symmetry breaking in graphs, Graph Theory Notes N.Y. 30 (1996) 6-7.

\bibitem{a_c} M. Albertson and K. Collins, Symmetry breaking in graphs, Electronic Journal of Combinatorics 3 (1996).

\bibitem{b_c} B. Bogstad and L. Cowen, The distinguishing number of the hypercube, Discrete Mathematics 283 (2004) 29-35.

\bibitem{mc_products} M. Chan, The distinguishing number of the direct product and wreath product action, Journal of Algebraic Combinatorics, to appear.

\bibitem{mc_maximum} M. Chan, The maximum distinguishing number of a group, Electronic Journal of Combinatorics, to appear.

\bibitem{cctcheng} C.C.T. Cheng, Three problems in graph labeling, Ph.D. Thesis, Department of Mathematical Sciences, Johns Hopkins University, 1999.

\bibitem{c_s} S. A. Choudum and V. Sunitha, Augmented cubes, Networks, 40 (2002), 71-84.

\bibitem{c_stechnical} S.A. Choudum and V. Sunitha, Automorphisms of augmented cubes, Technical report, Department of Mathematics, Indian Institute of Technology Madras, Chennai, March 2001.

\bibitem{archive} Open problems column, SIAM Activity Group Newsletter in Discrete Mathematics, Summer-Fall 1996, No. 23, archived at http://www.math.uiuc.edu/\~{}west/pcol/pcolink.html.

\bibitem{m_p} Z. Miller and M. Perkel, A stability theorem for the automorphism groups of powers of the $n$-cube,  Australasian Journal of Combinatorics 10 (1994), 17-28.

\bibitem{potanka} K. Potanka, Groups, graphs and symmetry breaking, Masters Thesis, Department of Mathematics, Virginia Polytechnic Institute, 1998.

\bibitem{r_s} A. Russell and R. Sundaram, A note on the asymptotics and computational complexity of graph distinguishability, Electronic Journal of Combinatorics 5 (1998).

\bibitem{jt} J. Tymoczko, Distinguishing numbers for graphs and groups, Electronic Journal of Combinatorics 11 (1) (2004).

\end{thebibliography}
